\documentclass[11pt,english]{smfart}

\usepackage[utf8]{inputenc}
\usepackage[T1]{fontenc}
\usepackage[english,french]{babel}

\usepackage{amsmath,amssymb,amscd,amsbsy,bm,bbm,mathrsfs,mathtools}
\usepackage{euscript}
\usepackage{faktor}
\usepackage{scalerel}
\usepackage{manfnt}
\usepackage{graphicx}
\usepackage{tikz-cd}

\usepackage[text={150mm,251mm},centering,marginparwidth=75pt]{geometry}
\usepackage[dvipsnames]{xcolor}
\usepackage{colortbl}
\usepackage{color}

\usepackage{enumitem}
\setlist[itemize,enumerate]{wide=0pt, labelwidth=2em, labelsep*=0em, itemindent=0pt, leftmargin=\dimexpr\labelwidth + \labelsep\relax, noitemsep, topsep=1ex}

\usepackage[colorlinks,
  linkcolor=BrickRed,
  citecolor=Green,
  urlcolor=Cerulean,
  hypertexnames=true]{hyperref}
\usepackage{smfhyperref}
\usepackage[hyperpageref]{backref}
\renewcommand{\backref}[1]{$\uparrow$~#1}

\usepackage{cleveref}
\usepackage{nameref}

\usepackage{datetime}
\usepackage{comment}
\usepackage{cite}
\usepackage{thmtools}
\usepackage[all]{xy}
\usepackage{marginnote}
\usepackage{interval}
\usepackage{letltxmacro}
\usepackage{xspace}

\makeatletter
\newcommand{\crefnames}[3]{%
  \@for\next:=#1\do{%
    \expandafter\crefname\expandafter{\next}{#2}{#3}%
  }%
}
\makeatother

\renewcommand{\thethmx}{\Alph{thmx}}
\newtheorem{theorem}{Theorem}[section]

\theoremstyle{definition}

\theoremstyle{remark}

\numberwithin{equation}{section}
\newlist{thmlist}{enumerate}{1}
\setlist[thmlist]{font=\normalfont, label=(\roman*), ref=\thetheorem.(\roman{thmlisti})}

\newlist{thmenum}{enumerate}{1}
\setlist[thmenum]{font=\normalfont, label=(\roman*), ref=\thethmx.(\roman{thmenumi})}

\newlist{corlist}{enumerate}{1}
\setlist[corlist]{font=\normalfont, label=(\roman*), ref=\thecorx.(\roman{corlisti})}

\addtotheorempostheadhook[theorem]{\crefalias{thmlisti}{thm}}
\addtotheorempostheadhook[lemma]{\crefalias{thmlisti}{lem}}
\addtotheorempostheadhook[corollary]{\crefalias{thmlisti}{cor}}
\addtotheorempostheadhook[proposition]{\crefalias{thmlisti}{proposition}}
\addtotheorempostheadhook[definition]{\crefalias{thmlisti}{dfn}}
\addtotheorempostheadhook[remark]{\crefalias{thmlisti}{rem}}
\addtotheorempostheadhook[assumption]{\crefalias{thmlisti}{assumption}}

\crefname{lemma}{Lemma}{Lemmas}
\crefname{conjecture}{Conjecture}{Conjectures}
\crefname{theorem}{Theorem}{Theorems}
\crefname{proposition}{Proposition}{Propositions}
\crefname{definition}{Definition}{Definitions}
\crefname{remark}{Remark}{Remarks}
\crefname{corollary}{Corollary}{Corollaries}
\crefname{corx}{Corollary}{Corollaries}
\crefname{claim}{Claim}{Claims}
\crefname{assumption}{Assumption}{Assumptions}
\crefname{thmx}{Theorem}{Theorems}
\crefname{main}{Main Theorem}{Main Theorems}

\newcommand{\C}{\mathbb{C}}
\newcommand{\R}{\mathbb{R}}

\newcommand{\fkk}{\mathfrak{k}}
\newcommand{\fka}{\mathfrak{a}}
\newcommand{\fkp}{\mathfrak{p}}
\newcommand{\norm}[1]{\left\lVert#1\right\rVert}

\renewcommand\thefootnote{\textcolor{MidnightBlue}{\arabic{footnote}}}

\begin{document}

\title[]{A unified calculation for Gromov norm of K\"ahler class of bounded symmetric domains}

\author{Yuan Liu}
\email{lexliu@hku.hk}
\address{Department of Mathematics, The University of Hong Kong, Pokfulam Road, Hong Kong}
\urladdr{https://sites.google.com/view/yuan-lius-website}

\date{\today}

\begin{abstract}
We provide a unified way to calculate the Gromov norm of the K\"ahler class of all (compact manifolds uniformized by) bounded symmetric domains. This was done for three classical domains by Domin and Toledo and for the general case by Clerc and \O rsted. Here, the calculation is much simplified by a combination of the ideas in Domin-Toledo and a work of Toledo, with the help of the Polydisc Theorem. The equality is achieved if and only if the triangle is ideal with three vertices on the Shilov boundary.
% Your abstract goes here
\end{abstract}

\maketitle

\begingroup
\renewcommand{\thefootnote}{}
\footnotetext{%
\noindent
\parbox{\textwidth}{%
\textit{2020 Mathematics Subject Classification.} 
Primary 53C35, 55M25. Secondary 53C55.

\textit{Keywords.} Gromov norm; bounded symmetric domain; K\"ahler class.
}%
}
\endgroup

\section{Introduction}
Let $X$ be an $n$-dimensional bounded symmetric domain of rank $r$, or equivalently, a Hermitian symmetric space of noncompact type. Let $\omega$ be the K\"ahler form corresponding to the Bergman metric $g_0$, with the metric normalized such that the holomorphic sectional curvature $\kappa$ satisfies $-1\le \kappa\le -\frac{1}{r}$. The calculation of the Gromov norm of $\omega$ was done in \cite{DomicToledo87} for domains of type I, II, III, and in \cite{ClercOrsted} for the general case. The result is
\begin{equation}\label{eq: main}
\|\omega\|_\infty= r\pi
\end{equation}
Since this value is additive with respect to the product structure of $X$, from now on, we assume that $X$ is \textit{irreducible}. Express $X$ as a coset space of Lie groups: $X=G_0/K$ where $G_0$ is the largest connected group of biholomorphic automorphisms of $X$ and $K$ is the isotropy group of a fixed point $o=eK\in X$.

According to \cite[Section 1]{DomicToledo87}, the calculation of the Gromov norm can be reduced to the estimation of the integral  
$$\int_T \omega$$
where $T$ is a geodesic triangle in $X$. We will start with this reduction and provide a unified way of estimating this integral. Assume that the three vertices of $T$ are $P, Q$, and $R$, and write it as $T=T(P, Q, R)$. We break the argument into the following steps:
\begin{enumerate}[label=(\theenumi)]
    \item Apply an action of $g\in G_0$ such that $g.P=o$, due to transitivity of group action.
    \item Apply the action of $k\in K$ to move $Q$ into a fixed polydisc $\Delta^r$, keeping the same notation for $Q$ and $R$.
    \item Take orthogonal projection $\pi(R)$ of $R$ into $\Delta^r$ along geodesics, and show that $\int_{T(o,Q,R)}\omega=\int_{T(o,Q,\pi(R))}\omega$.
    \item The question is now reduced to $\Delta^r$, and the estimation \cref{eq: main} is deduced from the calculation in $\Delta$ and the additive property of this invariant.
\end{enumerate}

The idea of projection used in step 3 is from another paper of Toledo \cite{Toledo89}, where he dealt with the rank 1 case, i.e., $X$ is the complex unit ball. Besides, the \textit{special} potential $\varrho_o$ of the Bergman metric centered at $o$ was introduced in \cite[p. 427--428]{DomicToledo87} (or see \cite[Section 2]{ClercOrsted}) for calculation, which will reduce the integral over a triangle $T(o, Q, R)$ to $\int_{\overline{QR}} d^{\C} \varrho_o$, with $\overline{QR}$ the geodesic segment connecting $Q$ and $R$. We will not repeat the argument, while the terminology and notation will be explained in \cref{subsec: old steps}. Since we also need the property of special potential functions in step 2 here, we introduce it before carrying out the four steps in \cref{subsec: new steps}.

This paper is organized as follows. In \cref{sec: reduction}, we give the reduction of Domin-Toledo from the Gromov norm to the integral on geodesic triangles. In \cref{sec: estimation of the integral}, we introduce the necessary background on bounded symmetric domains and estimate the integral. This includes the Polydisc Theorem in \cref{subsec: background of bsd}, the construction of special potential functions in \cref{subsec: old steps},
and the above four steps in \cref{subsec: new steps}. The condition of equality is stated in the end using the Gauss-Bonnet formula.

\section{Reduction to the integral on triangles}\label{sec: reduction}

In this section, we give the definition of the Gromov norm introduced in \cite{Gromov82} and the reduction of this norm done by Domin-Toledo in \cite[Section 1]{DomicToledo87}.

Let $X$ be a topological space and $c$ a singular cochain on $X$, the sup norm of $c$ is defined as
$$\norm{c}_{\infty}=\sup\{|c(\sigma)|\colon \sigma \text{\ a singular simplex in\ }X \}$$ 
If $\alpha\in H^*(X;\R)$, the sup norm of $\alpha$ is defined as
$$\norm{\alpha}_{\infty}=\inf\{\norm{c}_\infty\colon c \text{\ any singular cocycle representing\ } \alpha\}$$
which defines a pseudo-norm on cohomology classes.

In this note, we do the calculation for $X$, a bounded symmetric domain, and $[\omega]\in H^2(X;\R)$ the cohomology class of its canonical K\"ahler form $\omega$. We should evaluate any 2-form $\alpha$ representing $[\omega]$ on any 2-simplex $\sigma$, say $\int_\sigma \alpha$. First, notice that the integral is independent of the choice of $\alpha$, due to Green's Theorem. From now on, we work with $\omega$. Next, since $X$ is a complete manifold with strictly negative sectional curvature, we only need to consider for $\sigma$ a geodesic triangle, i.e., a (solid) triangle $T$ whose three edges are geodesics (see \cite[p. 426, (3)]{DomicToledo87}). Assume the three vertices of $T$ are fixed, say $P, Q, R$, then the three sides of triangle $T$ are uniquely determined due to the curvature condition. Since $\omega$ is closed, by Stokes' Theorem, the integral $\int_T \omega$ is independent of the choice of 2-cells filling up this triangle, say $T$. By the above argument, we henceforth focus on
$\int_T\omega$ with $T$ any geodesic triangle with three vertices labeled as $P, Q$ and $R$.

The above argument works verbatim for $X$ a compact K\"ahler manifold uniformized by a Hermitian symmetric space of noncompact type. This is due to $\omega$ is invariant under the biholomorphic action of the Hermitian symmetric space, as stated in \cite{DomicToledo87}. 

\section{Estimation of the integral}\label{sec: estimation of the integral}
\subsection{Background on the bounded symmetric domains}\label{subsec: background of bsd}

Let $X=G_0/K$ be an irreducible bounded symmetric domain, with $o=e.K\in X$ a fixed point. The property of bounded symmetric domains we will apply is the Polydisc Theorem (see \cite[p. 280]{Wolf}). The readers may consult \cite{Wolf} and \cite{Mok89} for more information on the fine structures of bounded symmetric domains. 

Writing $\mathfrak{g}_0$ and $\mathfrak{k}$ as the Lie algebras of $G_0$ and $K$ respectively, we have the Cartan decomposition $\mathfrak{g}_0=\mathfrak{k}+\mathfrak{p}$. Let $\fka$ be a Cartan subspace, i.e., a maximal abelian subspace of $\fkp$, and write $A=\exp(\fka)$. By the standard symmetric space theory, we have $G_0=KAK$, and this decomposition is consistent with the following Polydisc Theorem.

\begin{theorem}[Polydisc Theorem]\label{thm: polydisc theorem}
With the notation as above, there exists a totally geodesic complex submanifold $D$ of $X$ such that $(D,g_0|_D)$ is isometric
to a Poincaré polydisc $\Delta^r$ (equipped with product of the Poinca\'e metric) and $$X=\bigcup_{k\in K}k.D$$
\end{theorem}

The consistency of decomposition means that $X=G.o=KAK.o=KA.o$ and $D=A.o$. Sometimes we abuse the notation a little bit and write $D=\Delta^r$.

\subsection{Construction of special coordinates}\label{subsec: old steps}

Since $G_0$ acts transitively on $X$, we can take $P=o$. Note that $G_0$ acts isometrically on $X$ and the estimation of the integral $\int_T\omega$ for $T$ any geodesic triangle is not affected by $G_0$ action. From now on, we take $P=o$.

Next, we give the construction of special potential functions in \cite{DomicToledo87} (or refer to \cite[Section 2]{ClercOrsted}). There exists a unique potential function $\varrho_o$ for the Bergman metric, i.e., $dd^\C\varrho_o=\omega$ with the properties follows:
\begin{enumerate}[label=(\roman*)]
    \item $\varrho_o(o)=0$.
    \item $\varrho_o$ is invariant under $K$-action.
    \item $d^\C\varrho_o$ vanishes on the tangent vector to any geodesic through $o$.
\end{enumerate}
here $d^\C=\sqrt{-1}(\bar{\partial}-\partial)$. We call this $\varrho_o$ the \textit{special potential function} of the Bergman metric. The first two properties are rather obvious, and the readers may refer to \cite{DomicToledo87,ClercOrsted} for (iii). By (iii), we have
$$\int_{T(o,Q,R)}\omega=\int_{\partial T(o,Q,R)}d^\C \varrho_o=\int_{\overline{QR}} d^\C \varrho_o$$
where $\overline{QR}$ means the geodesic segment from $Q$ to $R$.

For any different point $A\in X$, let $g\in G_0$ such that $A=g.o$. Then define 
\begin{equation}\label{eq: def rho A}
\varrho_A(z)=\varrho_o(g^{-1}.z)   
\end{equation}
which satisfies $dd^\C\varrho_A=\omega$ and
\begin{enumerate}[label=(\roman*)]\label{item total}
    \item $\varrho_A(A)=0$.
    \item \label{item:a} $\varrho_A$ is invariant under $\mathrm{Ad}(g)K$-action.
    \item $d^\C\varrho_A$ vanishes on the tangent vector to any geodesic through $A$.
\end{enumerate}

This is the special potential function at $A$ for the Bergman metric. Take $A=Q$ and $u=\varrho_o-\varrho_Q$, then $u$ is pluriharmonic and we can write $u=\mathrm{Re}(H)$ for $H$ a holomorphic function on $X$, due to $X$ is contractible. Take $v=\mathrm{Im}(H)$, then 
\begin{align*}
\int_{T(o,Q,R)}\omega=&\int_{\overline{QR}} d^\C \varrho_o=\int_{\overline{QR}} d^\C \varrho_o-d^\C\varrho_Q\\
=&\int_{\overline{QR}}d^\C u=\int_{\overline{QR}}d v=v(Q)-v(R)  
\end{align*}
Thus, to prove \cref{eq: main}, we only need to show that 
\begin{equation}\label{eq: imaginary part}
|v(Q)-v(R)|\le r\pi
\end{equation}

\subsection{Estimation of the integral}\label{subsec: new steps}

Now let us go back to $\int_T \omega$, with
$T=T(o,Q,R)$ a geodesic triangle. By the \cref{thm: polydisc theorem}, we have $k\in K$ such that $k.Q\in \Delta^r$, where $\Delta^r$ is a specified polydisc containing $o$, fixed from now on. We keep the same notation $Q, R$ for the image of $Q, R$ under this action of $k$. Thanks to the property (ii) of $\varrho_o$ in \cref{subsec: old steps}, the function $\varrho_o$ is unchanged. Now take the special potential function $\varrho_Q$ for this new $Q$, and continue the process in \cref{subsec: old steps}.

The new triangle $T(o, Q, R)$ does not necessarily belong to $\Delta^r$, and we denote by $\pi: X\to \Delta^r$ the orthogonal projection onto $\Delta^r$ along geodesics. This operation was considered by Toledo in \cite[p. 128]{Toledo89}, when $r=1$, i.e., $X=B^n$ the unit complex ball in $\C^n$. Now, we follow his idea to show that 
\begin{equation}\label{eq: projection equal}
\int_{T(o,Q,R)}\omega=\int_{T(o,Q,\pi(R))}\omega
\end{equation}

Denote the tetrahedron formed by $o,Q,R,\pi(R)$ as $V$, then by Stokes' theorem, we have 
$$\int_{\partial V}\omega =\int_V d\omega=0$$ 
Thus, to prove \cref{eq: projection equal}, we simply need to show the vanishing of the integral of $\omega$ on any of the two faces containing the geodesic segment $\overline{R\pi(R)}$.

For $\overline{oQ}\in \Delta^r\subseteq\C^r$ equipped with the canonical complex structure, consider the projection from $X$ onto $\Delta^r$ along geodesics. The geodesic connecting $R$ and $\pi(R)$ is orthogonal not only to $\overline{o\pi(R)}$, but also to the whole complex plane $\C^r$ containing $\overline{oQ}$. The restriction of $\omega$ to $T(o, R, \pi(R))$ is identically 0, and we have the following two ways of seeing this. The face $T(o, R,\pi(R))$ is contained in a totally real submanifold of $X$, and thus the restriction of $\omega$ there is identically 0. This is because $\C^r=\exp(\fka)$ and the orthogonal complement of the tangent space of $\C^r$ at a given point $\pi(R)$ is identified with $\fkk$. Or in another way, we can show the vanishing of $\omega$ along $\overline{R\pi(R)}$ using the property (ii) of $\varrho_o$ in \cref{subsec: old steps} and the fact $\int_{T(o, R, \pi(R))}\omega=\int_{\overline{R\pi(R)}}d^\C\varrho_o$. The same holds true for the other face $T(Q, R,\pi(R))$, and \cref{eq: projection equal} is proved.

Now we are reduced to the situation that the triangle $T(o, Q, R)$ is contained in $\Delta^r$, which we assume from now on. Recall that in \cref{subsec: old steps}, the estimation of the integral over a triangle is reduced to \cref{eq: imaginary part} and is additive with respect to the product structure of the space. 
The $\Delta^r\subset X$ is totally geodesic, being isometric to $(\Delta,g_0)^r$ where $g_0$ is the Poincaré metric, normalized such that the curvature of each factor is $-1$. In conclusion, we only need to show that 
$$\int_{T(o, Q, R)} \omega\le \pi$$
for $T(o, Q, R)$ a geodesic triangle in $\Delta$. The rest is direct calculation and is the special case of $p=q=1$ for the domain of type I in \cite[Section 2]{DomicToledo87} (see also \cite{Toledo89}). We include this here for completeness and provide another way of seeing this with the Gauss-Bonnet formula later.

Take $z$ as the coordinate of the unit disc $\Delta$ and two points $Q, R\in \Delta$ with coordinates $z_1$ and $z_2$, respectively. Now the special potential function at the origin $o$ is 
$$\varrho_o(z)=-\log(1-|z|^2)$$
The M\"obius transformation $g^{-1}$ moving $Q$ to $o$ is
$$g^{-1}.z=\frac{z-z_1}{-\bar{z}_1z+1}$$
By \cref{eq: def rho A}, we have 
$$\varrho_Q=\varrho_o\left(\frac{z-z_1}{-\bar{z}_1z+1}\right)=-\log(1-|z|^2)-\log(1-|z_1|^2)+\log(1-|\bar{z}_1z|^2)$$
Thus $u=\varrho_o-\varrho_Q=\log\left(\frac{1-|z_1|^2}{|1-\bar{z}_1z|^2}\right)$ and its harmonic conjugation is 
$$v=-2\arg (1-\bar{z}_1z)$$
Writing $\bar{z}_1{z}_2=\rho e^{i\theta}$ with $0<\rho<1$, we have 
$$|v(z_2)-v(z_1)|=2|\arg(1-\rho e^{i\theta})|<\pi$$
We have now proved that $\norm{\omega}_\infty\le r\pi$. To approach the upper bound $\pi$, roughly speaking, we need to take $\rho\to 1^-$, and then $\theta\to 0$ (or $2\pi$), which means that the geodesic triangle $T$ is ideal, i.e., all three vertices are on the boundary $\partial \Delta$. The readers may refer to \cite{Toledo89} for more details of this argument of equality. 

We can also see this using the Gauss-Bonnet formula as follows. The (normalized holomorphic) sectional curvature of $\Delta$ equipped with the Poincaré metric is $-1$. For $T$ the solid triangle contained in $\Delta$, we have that
$$\left|\int_T\omega\right|=\mathrm{Area}(T) =\pi-(\alpha+\beta+\gamma)\le \pi$$
where $\alpha,\beta,\gamma$ are the interior angles of $T$. It is $\pi$ if and only if the triangle is ideal. This gives the condition of equality in \cref{eq: main}; the triangle must have three vertices on the boundary of each component of $\Delta^r$. With the notation in the 4 steps, what we have shown is that $T(o, Q, \pi(R))$ is ideal with all three vertices on the Shilov boundary of $X$. By Hermann Convexity Theorem (\cite[p. 286]{Wolf}), the vertex $R$ lies also on the boundary.

\section*{Acknowledgments}
The author is partially supported by the National Natural Science Foundation of China (Grant No. 12501105). The author thanks Professor N. Mok for many useful conversations and encouragements during this work.

\bibliographystyle{ssmfalpha}
\bibliography{reference}

\end{document}